\documentclass[11pt,a4paper,twoside]{article}
\usepackage{amsmath,amsfonts,amssymb,amsthm,graphics,hyperref}
\usepackage[mathscr]{euscript}

\newtheorem{Theorem}[equation]{Theorem}
\newtheorem{Corollary}[equation]{Corollary}
\newtheorem{Proposition}[equation]{Proposition}
\newtheorem{Lemma}[equation]{Lemma}

\newtheorem{Remarks}[equation]{Remarks}

\def\Section#1{\section{#1}\setcounter{equation}{0}}

\font\smc=cmcsc10

\def\wt{\widetilde}
\def\ov#1{\overline{#1}}

\def\la{\left\langle}
\def\ra{\right\rangle}

\def\spcheck{^\vee}

\def\bdots{\mathinner{\mkern1mu\raise2pt\hbox{.}\mkern2mu\raise5pt\hbox{.}
   \mkern2mu\raise8pt\hbox{.}\mkern1mu}}

\def\today{\number\day\space
 \ifcase\month\or
	January\or February\or March\or April\or May\or June\or
	July\or August\or September\or October\or November\or December\fi
 \space\number\year}

\def\oB{\ov{B}{}}
\def\oG{\ov{G}{}}
\def\oH{\ov{H}{}}
\def\oJ{\ov{J}{}}
\def\oK{\ov{K}{}}
\def\oL{\ov{L}{}}

\def\oP{\ov{P}{}}
\def\oU{\ov{U}{}}
\def\ol{\ov{\lambda}{}}

\def\pif{{\varpi_F}}
\def\pie{{\varpi_E}}

\def\vph{\varphi}
\def\of{\mathfrak o_F}
\def\pf{\mathfrak p_F}
\def\oe{\mathfrak o_E}

\def\fa{\mathfrak a}

\def\fs{\mathfrak s}
\def\o{\mathfrak o}
\def\p{\mathfrak p}
\def\A{\mathfrak A}
\def\B{\mathfrak B}
\def\H{\mathfrak H}
\def\J{\mathfrak J}

\def\P{\mathfrak P}
\def\Q{\mathfrak Q}

\DeclareMathOperator{\Aut}{Aut}
\DeclareMathOperator{\End}{End}
\DeclareMathOperator{\Hom}{Hom}

\DeclareMathOperator{\Ind}{Ind}
\DeclareMathOperator{\Res}{Res}
\DeclareMathOperator{\Ker}{Ker}

\DeclareMathOperator{\Lie}{Lie}
\DeclareMathOperator{\diag}{diag}

\def\det{\hbox{\rm det}}

\def\cB{\mathscr B}
\def\cC{\mathcal C}

\def\cG{\mathcal G}
\def\cH{\mathcal H}
\def\cL{\mathcal L}

\def\sm{{\mathsf m}}
\def\sM{{\mathsf M}}
\def\sW{{\mathsf W}}

\def\bC{\mathbb C}
\def\bM{\mathbb M}
\def\bN{\mathbb N}

\def\bZ{\mathbb Z}

\def\bI{\boldsymbol 1}
\def\bg{\boldsymbol g}
\def\bk{\boldsymbol k}
\def\bG{\boldsymbol G}

\def\bB{\boldsymbol B}
\def\bL{\boldsymbol L}
\def\bN{\boldsymbol N}
\def\bP{\boldsymbol P}
\def\bT{\boldsymbol T}
\def\bU{\boldsymbol U}
\def\bSp{\boldsymbol S\boldsymbol p}
\def\bSO{\boldsymbol S\boldsymbol O}

\def\bRes{\mathbf R\mathbf e\mathbf s}
\def\bv{\mathbf v}
\def\bw{\mathbf w}

\def\boG{\mathbf G}
\def\boGLMF{\mathbf G\mathbf L_{\mathbf M}\mathbf(\mathbf F\mathbf)}

\def\a{\alpha}
\def\b{\beta}

\def\d{\delta}
\def\e{\varepsilon}
\def\z{\zeta}
\def\th{\theta}
\def\k{\kappa}
\def\l{\lambda}
\def\s{\sigma}
\def\t{\tau}
\def\om{\omega}

\def\L{\Lambda}

\def\Om{\Omega}

\begin{document}

\title{Covers for self-dual supercuspidal representations of the\\ Siegel
Levi subgroup of classical $p$-adic groups}
\author{\setcounter{footnote}{2}%
David Goldberg\thanks{The first author was supported by NSF grant DMS0400958.}, 
Philip Kutzko\thanks{The second author was supported by NSF grant DMS0101451.} \ and 
\setcounter{footnote}{1}%
Shaun Stevens\thanks{The third author was supported by grant NAL/00779/G from the Nuffield Foundation.}} 
\date{\today}
\maketitle
\begin{abstract}
We study components of the Bernstein category for a $p$-adic classical 
group (with $p$ odd) with inertial support a self-dual positive level supercuspidal representation of 
a Siegel Levi subgroup.  
More precisely, we use the method of covers to construct a Bushnell-Kutzko type for such a component.
A detailed knowledge of the Hecke algebra of the type should have number-theoretic implications.
\end{abstract}


\section*{Introduction}

Let $F_0$ be a locally compact nonarchimedean local field, and $G$ the $F_0$-points of a
connected reductive quasi-split algebraic group $\bG$ defined over
$F_0$. The smooth representation theory of $G$ plays a role in
automorphic forms, and is therefore of deep interest to number
theorists. A theorem of Bernstein shows that, in the terminology of~\cite{BK1}, it is enough to
understand the subcategories $\mathfrak R^{\fs}(G)$ of smooth
representations with supercuspidal support in a given inertial class
$\fs$. The class $\fs=[L,\pi]_G$ is determined by a Levi subgroup
$L\subset G$ and an irreducible supercuspidal representation $\pi$ of
$L$. The theory of types and covers~\cite{BK1}
gives one a method for trying to understand the subcategories $\mathfrak
R^{\fs}(G)$. In particular, if one knows supercuspidal types for all
Levi subgroups of $G$, and
how to construct $G$-covers of all such types, then one has a
complete set of types for $G$. Here, we study the case where $F$
has odd residual characteristic,
$G$ is a classical group, and $\fs=[L,\pi]_G$ with $L$ a Siegel Levi
subgroup and $\pi$ an irreducible supercuspidal representation fixed
by the non-trivial Weyl group element. While this work gives a
complete set of $G$-covers (and hence types) in this case, the theory
cannot be considered complete without knowledge of the associated
Hecke algebras. As shown in~\cite{KM}, a precise understanding of 
such Hecke algebras gives
information on reducibility of induced representations and therefore
determines poles of certain Langlands $L$-functions. These
considerations will be the subject of future work. 

\medskip

Let $\bG$ be any of the groups $\bSp_{2M},\, \bSO_{2M},\, \bSO_{2M+1},\,
\bU_n$ or $\bSO_n^*$, i.e., a symplectic, quasi-split special
orthogonal group, or quasi-split unitary group. (Here $\bSO_n^*$ stands
for a quasi-split but non-split special orthogonal group.) Let
$\bB=\bT\bN$ be a Borel subgroup with unipotent radical $\bN$ and
maximal torus $\bT$. There is a choice of standard maximal
parabolic subgroup $\bP=\bL\bU$ containing $\bB$ whose Levi component
$\bL$ is isomorphic to either $\bG\bL_M\times\bG_0$, or $\bRes_{F/F_0}(\bG\bL_M)\times \bG_0$; here
$\bG_0$ is an anisotropic group and $F/F_0$ is the quadratic extension
defining the unitary group. We study the theory of types for the
inertial class $\fs=[L,\pi]_{G}$, 
where $\pi$ is fixed by the non-trivial element $w_0$ of the Weyl
group $N_{\bG}(\bL)/Z_{\bG}(\bL)$. The work of Bushnell and Kutzko~\cite{BK} gives one a type $(J_L,\l_L)$ for the inertial class
$\fs_L=[L,\pi]_L$. 
The purpose of this paper is to construct a cover $(J,\l)$ of
$(J_L,\l_L)$ for each group $\bG$.
(This was done by Blondel~\cite{B2} when $\bG=\bSp_{2M}$; her methods are somewhat different.)
Therefore, by~\cite{BK1},the pair
$(J,\l)$ is a type for $\fs$. One then knows that the category $\mathfrak
R^{\fs}(G)$ is isomorphic to the category $\cH(G,\l)\text{--Mod}$ of
unital (left) modules over the $\l$-spherical Hecke algebra. 

This isomorphism of categories has implications both for the location of poles of certain $L$-functions, and for the classification of local galois representations. One example of this may be seen in recent work of Kutzko and Morris~\cite{KM}, where $\bG$ is one of the groups $\bSp_{2M}$, $\bSO_{2M}$ or $\bSO_{2M+1}$, and $\pi$ is a level zero self-dual supercuspidal representation. Specific information about the Hecke algebra $\cH(G,\l)$ is employed to give a purely local proof of Shahidi's theorem on the reducibility of parabolic induction~\cite{Sh}, which can also be described in terms of the poles of $L$-functions. 
On the other hand, recent work of Henniart~\cite{He2}, building on work of Harris and Taylor~\cite{HT} and Henniart~\cite{He1}, relates this to the classification of local galois representations by their images (e.g. symplectic, orthogonal). In the case of the level zero representations in~\cite{KM}, 
the classification of galois representations obtained via the Hecke algebra isomorphism was known already. However, in the general situation of this paper, the corresponding classification of galois representations is not fully understood and has been sought for some time. We hope to address this problem in a sequel.

\medskip

We now give a summary of the contents of the paper. 
In~\S\ref{S1.1}, we review the theory of types for self-dual representations, via results of
Blondel, which is crucial to the construction (see \cite[Proposition 2.2]{B2}). 
In~\S\ref{S1.classical}, we fix a self-dual representation $\pi$ of $L$, and use these results to define a well-adapted inner product on an $F$-vector space $V$, which defines the group $G$.
In~\S\ref{S1.3}, we show that there is a choice of \emph{simple stratum} (see~\cite{BK}), defining the type in $\pi$, which is particularly well suited. In particular, it gives rise to a \emph{skew semisimple} stratum (see~\cite{S4}) in $\End_F(V)$, which underpins the construction of $G$-cover.

We construct the $G$-cover $(J,\l)$ in~\S\ref{S2}, following the recipe of~\cite[\S7]{BK}. The first stages, when we are working with pro-$p$ groups, are performed first in $\Aut_F(V)$ (see \S\ref{S2.1}), and then transferred to $G$ using the Glauberman correspondence (see \S\ref{S2.2}).
To verify that $(J,\l)$ is indeed a $G$-cover, we construct an invertible element $f$ in the Hecke algebra $\cH(G,\l)$, which is  supported on a strongly $(P,J)$-positive element in the
centre of $L$ (see~\cite{BK1}). To do this, we use
an argument which harks back to the work of Borel~\cite{Bo}: we find two invertible elements of the Hecke algebra $\cH(G,\l)$, which each have support in a compact subgroup of $G$, and whose convolution is supported on a strongly positive element of $L$ (see Lemmas~\ref{2.9},~\ref{2.10}); a suitable power of this is the required element $f$.

Finally, in~\S\ref{S2.3}, we prove a result on the Hecke algebra $\cH(G,\l)$, following the techniques and philosophy of~\cite[\S5]{BK} (see also~\cite{Ro}). 
We show that the subalgebra of elements with support in a fixed maximal compact subgroup 
is isomorphic to
an algebra of the form $\cH(\cG',\rho')$, where $\cG'$ is a (possibly
disconnected) classical group over a finite field, 
and $\rho'$ is an irreducible cuspidal representation of the Siegel
Levi subgroup of $\cG'$. In
many cases, the calculations in~\cite{KM} will give the parameters of
these Hecke algebras; in general, the calculations will be
similar. This will simplify the computation of Hecke algebra $\cH(G,\l)$;
these specific calculations, and their implications to the classification of galois representations, are left to future work.




\Section{Preliminaries}\label{S1}

Let $F$ be a locally compact nonarchimedean local
field. Let $\mu$ be an automorphism of $F$ with $\mu^2=1$; we allow 
the possibility that $\mu$ is trivial. We set
$F_0=F^\mu$ to be the fixed points of $\mu$. Let $\o_F$ be the ring
of integers in $F$, and $\p_F$ the maximal ideal in $\o_F$. Denote by
$k_F$ the residual field $\o_F/\p_F$, and let $q_F=|k_F|$. We adopt
similar notation for $F_0$ and for any extension of $F_0$. 

Let $p$ denote the characteristic of $k_F$. 
We assume that $p$ is not $2$ throughout.

For $r$ a real number, we write $\lfloor r\rfloor$ greatest integer 
less than or equal to $r$, and $\lceil r\rceil$ for the least integer 
greater than or equal to $r$.


\subsection{Self-dual representations of $\boGLMF$}\label{S1.1}

We begin by looking at the self-dual representations of $GL_M(F)$,
following Blondel~\cite[\S 2.2]{B2}.
Let $W$ be an $M$-dimensional $F$-vector space with
basis $\cB=\{\bw_1,\ldots,\bw_M\}$, equipped with the non-degenerate hermitian
form $\la\,,\,\ra_W$ given by $\la \bw_i,\bw_j\ra_W=\d_{i+j,M+1}$. Denote by
$\wt\e_W$ the adjoint involution on $A_W=\End_F(W)$ and by
$\e_W$ the involution $g\mapsto \wt\e_W(g^{-1})$ of $G_W=\Aut_F(W)\simeq
GL_M(F)$. 
Then, writing $a\in A_W$ with respect to the basis $\cB$, 
we have 
$$
\wt\e_W(a)=
{}^\dag a^\mu,
$$
where $a^\mu\in A_W$ is obtained by applying $\mu$ to the coefficients
of $a$, and ${}^\dag$ denotes transpose with respect to the off-diagonal.
For $\rho$ a representation of a subgroup $J$ of $G_W$, we
write $\rho^{\e_W}$ for the representation of $\e_W(J)$ given by
$\rho^{\e_W}(\e_W(j))=\rho(j)$, for $j\in J$.

Let $\pi$ be an irreducible supercuspidal representation of
$G_W$ such that $\pi\simeq\pi^{\e_W}$. Note that, if $F=F_0$, then, by
Gel$'$fand and Kazhdan~\cite{GK}, this is equivalent to 
$\pi$ being self-contragredient. Let
$(J_W,\l_W)$ be a maximal simple type in $G_W$ corresponding to the
inertial equivalence class $[G_W,\pi]_{G_W}$. We will need to use the
construction of $\l_W$ quite explicitly so we recall it briefly here
(see~\cite{BK} for more details).

It begins with a \emph{simple stratum} $[\A_W,n_W,0,\b]$, where $\A_W$
is a principal hereditary $\of$-order in $A_W$, with Jacobson radical
$\P_W$, and $n_W\in\mathbb N$ is such that $\b\in\P_W^{-n_W}\setminus
\P_W^{1-n_W}$. Further, $E=F[\b]$ is a field extension of $F$ and
$E^\times$ normalizes $\A_W$. We set $B_W$ to be the $A_W$-centralizer
of $E$ and put $\B_W=\A_W\cap B_W$, a maximal hereditary $\oe$-order
in $B_W$. We also fix a uniformizer $\pie$ of $E$.

From the stratum are defined certain subgroups $H^k_W=H^k(\b,\A_W)$
and $J^k_W=J^k(\b,\A_W)$, for $k\ge 0$, along with some sets
$\cC(\A_W,k,\b)$ of characters of $H_W^{k+1}$ called \emph{simple
characters} (see~\cite[\S3]{BK}). The construction of the type continues
with a simple character $\th_W\in\cC(\A_W,0,\b)$. There is then a
unique irreducible representation $\eta_W$ of $J^1_W$ which contains
$\th_W$. 

Now we take $\k_W$ to be a \emph{$\b$-extension} of $\eta_W$, that
is, one of a certain family of representations of $J_W=J^0_W$ which
restrict to $\eta_W$. We also recall the construction of $\k_W$ here,
from~\cite[\S5.1--2]{BK}. Let $\B^\sm_W$ be a minimal hereditary $\oe$-order
in $B_W$ contained in $\B_W$ and let $\A^\sm_W\subset\A_W$ be the unique
hereditary $\of$-order in $A_W$ which is normalized by $E^\times$ and
such that $\A^\sm_W\cap B_W=\B^\sm_W$. Then $[\A^\sm_W,n^\sm_W,0,\b]$ is also
a simple stratum, for some integer $n^\sm_W$, and we can define simple
characters associated to this stratum also. Moreover, there is a
canonical bijection $\t_{\A_W,\A^\sm_W,\b}:\cC(\A_W,0,\b)\to
\cC(\A^\sm_W,0,\b)$ (see~\cite[\S3.6]{BK}). Let $\th_W^\sm$ be the transfer of
$\th_W$ under this bijection. There is a unique irreducible
representation $\eta^\sm_W$ of $J^1(\b,\A^\sm_W)$ which contains
$\th^\sm_W$.

Now we form the group $\wt J^1_W=U^1(\B^\sm_W)J^1_W$. There is a
unique representation $\wt\eta_W$ of $\wt J^1_W$ such that
$\wt\eta_W|_{J^1_W}=\eta_W$ and $\wt\eta_W,\eta^\sm_W$ induce
equivalent irreducible representations of $U^1(\A_W^\sm)$. Then $\k_W$
is any representation of $J_W$ such that $\k_W|_{\wt J{}^1_W}=\wt\eta_W$. 
Finally, $J_W/J^1_W\cong U(\B_W)/U^1(\B_W)$ is isomorphic to
$GL_r(k_E)$, where $r=M/[E:F]$. Then there is a cuspidal
representation $\rho_W$ of $J_W/J^1_W$ such that
$$
\l_W=\k_W\otimes\rho_W.
$$

\medskip

By conjugating $(J_W,\l_W)$ if necessary, we may and do assume $\A_W$
is \emph{standard}; that is, with respect to our chosen basis $\cB$, it
consists of matrices with entries in $\of$ which are upper block
triangular modulo $\pf$. Note that this means that $\wt\e_W(\A_W)=\A_W$.

\begin{Proposition}[{cf.~\cite[2.2 Proposition]{B2}}]\label{1.2}
\begin{enumerate}
\item There exists $\s\in U(\A_W)$ such that $J_W$ is stable
under $\wt\s: g\to \s\, \e_W(g) \s^{-1}$ and $\l_W$ is equivalent to
$\l_W\circ\wt\s$.
\item Such an element $\s$ is unique up to left multiplication
by $J_W$. It satisfies: 
\begin{enumerate}
\item $\s\,\e_W(\s)\in J_W$ and $\pie^{-1}\s\, \e_W(\pie^{-1}\s) 
\in J_W$. 
\item The map $\wt\s$ stabilizes $H_W^1$ and $J_W^1$ and we
have $\th_W=\th_W\circ\wt\s$. 
\item The lattices $\J_W$, $\H^1_W$ are stable under
$X\mapsto \s\wt\e_W(X)\s^{-1}$. 
\end{enumerate}
\end{enumerate}
\end{Proposition}

The proof is identical to that of~\cite[2.2 Proposition]{B2}.


\subsection{Simple characters}\label{S1.2}

We record here, the following useful lemma, from~\cite[4.3 Lemma]{B2}.
Note that this uses very strongly the condition that $p\ne 2$.

\begin{Lemma}[{\cite[4.3~Lemma~1]{B2}}] \label{1.3}
Suppose $V''$ is any $F$-vector space and
$[\A'',n'',0,\b'']$ is
a simple stratum in $A''=\End_F(V'')$, where $\A''$ is a hereditary
$\of$-order in $V''$.
\begin{enumerate}
\item $[\A'',n'',0,\frac 12\b'']$ is a simple stratum in $A''$, with
$H^k(\frac 12\beta'',\A'')=H^k(\beta'',\A'')$, for each $k\ge 0$, and
similarly for $J^k$. 
\item For each $m\ge 0$, the map $\theta\mapsto\theta^2$ is a bijection
from $\cC(\A'',m,\frac 12\beta'')$ onto $\cC(\A'',m,\beta'')$ which is
compatible with the canonical bijections $\tau$ of~\cite[\S3.6]{BK}. 
We denote the inverse bijection by
$\theta\mapsto\th^{1/2}$.
\end{enumerate}
\end{Lemma}

We will write $\th_L=\theta_W^{1/2}\in
\cC(\A_W,0,\frac 12\beta)$ and $\th^\sm_L=(\theta^\sm_W)^{1/2}\in
\cC(\A^\sm_W,0,\frac 12\beta)$. We also let $\eta_L$ be the unique
irreducible representation of $J^1(\b,\A_W)$ which contains $\th_L$,
and likewise $\eta^\sm_L$.

Note also that we have $(\th_L\circ\wt\s)^2=\th_L^2\circ\wt\s = 
\th_W$. In particular, since 
the squaring map is a bijection, we see that $\th_L\circ\wt\s=\th_L$.


\subsection{Classical groups}\label{S1.classical}

Let $V_0$ be an $F$-vector space equipped with a nondegenerate
totally isotropic $\nu$-hermitian form $\la \,,\,\ra_0$, with $\nu
=\pm1$. Thus $\la \bv_0,\bw_0\ra_0=\nu\la \bw_0,\bv_0\ra_0^\mu$, for
all $\bv_0,\bw_0\in V_0$.  
We write $G_0=\Aut_F(V_0)$ and denote by $\oG_0^+$
the (anisotropic) group corresponding to the form
$\la\,,\,\ra_0$:
$$
\oG_0^+ := \{g_0\in G_0: \la g_0\bv_0,g_0\bw_0\ra_0=\la \bv_0,\bw_0\ra_0\hbox{
for all }\bv_0,\bw_0\in V_0\}.
$$
We put 
$
\oG_0 := \{g_0\in \oG_0^+: \det_{G_0/F}(g_0)=1\}.
$
We allow the possibility that $V_0=\{0\}$. 

Set $V=W\oplus V_0\oplus W$ and define a form $\la \,,\,\ra$ on $V$ by
$$
\la (\bv_1,\bv_0,\bv_2),(\bw_1,\bw_0,\bw_2)\ra = \la \bv_1,\s^{-1}\bw_2\ra_W+ \la \bv_0,\bw_0\ra_0 +\nu\la \s^{-1}\bv_2,\bw_1\ra_W,
$$
for $\bv_1,\bv_2,\bw_1,\bw_2\in W$ and $\bv_0,\bw_0\in V_0$, where
$\s$ is the element given by Proposition~\ref{1.2}.
Note that $\la\, ,\,\ra$ is now a nondegenerate $\nu$-hermitian form
in which the two copies of $W$ are (dual) maximal isotropic spaces,
and $V_0$ is the maximal anisotropic space. Let $N=\dim_F V$; then $N=2M+D$, where $D=\dim_F V_0$. 

We put $A=\End_F(V)\simeq\bM(N,F)$. We set
$G=A^\times= \Aut_F(V) \simeq GL_N(F)$ and 
$$
\oG^+ = \{g\in G: \la g\bv,g\bw\ra=\la \bv,\bw\ra\hbox{ for all }\bv,\bw\in V\},
$$
a unitary, symplectic or orthogonal group over $F_0$. We also put
$$
\oG := \{g\in \oG^+: \det_{G/F}(g)=1\}.
$$
More generally, for $H$ a subgroup of $G$, we will write $\oH^+$ for
the intersection $\oH^+=H\cap\oG^+$, and $\oH$ for the intersection
$\oH=H\cap\oG$. 

We denote by $\wt\e$ the adjoint involution on $A$ determined
by the form; that is, for $a\in A$, $\wt\e(a)$ is the unique
element of $A$ such that
$$
\la av,w \ra = \la v,\wt\e(a)w\ra,\qquad\hbox{ for all }v,w,\in V,
$$
We have an involution $\e$ on $G$ given by
$\e(g)=\wt\e(g^{-1})$, for $g\in G$, so that $\oG^+=G^\e$.

Note that $\tilde\e$ induces an involution on $A_W$ via the two embeddings
of $W$ into the first (respectively last) factor of $V=W\oplus V_0\oplus
W$, and we also denote this involution by $\tilde\e$. For any
$\bv,\bw\in W$ and $a\in A_W$, we have 
$\la a(\bv,0,0),(0,0,\bw)\ra = \la
(\bv,0,0),\tilde\e(a)(0,0,\bw)\ra$. So, by definition of the form
$\la\,,\,\ra$,  
$$
\la av,\s^{-1}w\ra_W\ =\ \la v,\s^{-1}\tilde\e(a)w\ra_W,
\text{ that is } 
\la v,\tilde\e(a)_W\s^{-1}w\ra_W\ =\ \la v,\s^{-1}\tilde\e(a)w\ra_W.
$$
Thus, for any $a\in A_W$, we have $\tilde\e(a)=\s\tilde\e_W(a)\s^{-1}$.

Given a representation $\rho$ of a subgroup $J$ of $G_W$, we denote by
$\rho^\e$ the representation of $\e(J)$ given by
$\rho^\e(\e(j))=\rho(j)$, for $j\in J$.

For $L$ an $\of$-lattice in $V$, we define its dual lattice to be
$$
L^\#=\{v\in V:\la v,L\ra\subset\pf\}.
$$
We note that, since $\oG^+_0$ stabilizes lattices $L_0\supset
L_0^\#\supset \pf L_0$ in $V_0$ (see~\cite[\S1.8]{M1}), it is contained in a
maximal $\e$-stable compact open subgroup of $G_0$ (namely, the
$G_0$-stabilizer of these lattices).

Let $P$ be the parabolic subgroup of $G$ stabilizing the self-dual
flag 
$$
\{0\}\subsetneq W\oplus\{0\}\oplus\{0\} \subseteq
\left(W\oplus\{0\}\oplus\{0\}\right)^\perp=
W\oplus V_0\oplus \{0\}\subsetneq V,
$$
with unipotent radical $U$. Let $L$ denote the Levi component
of $P$ which stabilizes each copy of $W$ along with $V_0$. So $L\cong
GL_M(F)\times GL_D(F)\times GL_M(F)$. Let $P_-$ denote the
opposite parabolic subgroup, $P_-=LU_-$. We also put $A_{L}=A_W\oplus
A_0\oplus A_W = \Lie(L)$.

We set $\oP=P\cap\oG$, the Siegel parabolic subgroup of $\oG$,
with unipotent radical $\oU=U\cap\oG$, and $\oL=L\cap\oG$, a Levi
component of $\oP$. Then $\oL\cong GL_M(F)\times\oG_0$; in block matrix
form, we identify $G_W\times\oG_0$ with $\oL$ via the isomorphism 
$$
i(g,g_0)= \begin{pmatrix}
g&0&0\\ 0&g_0&0\\ 0&0&\e(g) 
\end{pmatrix},
\qquad\hbox{for }g\in G_W,\ g_0\in\oG_0.
$$
If $\rho$ is a representation of a subgroup $J_W$ of $G_W$, and $\oJ_0$
is a subgroup of $\oG_0$, then we denote by $i(\rho)$ the
representation of $i(J_W\times\oJ_0)$ given by $i(\rho)(i(j,j_0))=\rho(j)$, for
$j\in J_W$, $j_0\in\oJ_0$. If the group $\oJ_0$ is not specified then we
take it to be the whole of $\oG_0$. 

More generally, if $\oJ$ is a subgroup of
$\oG$ such that $\oJ\cap\oL=i(J_W\times\oJ_0)$ and $\oJ$ has an Iwahori
decomposition with respect to $(\oL,\oP)$, then we denote by
$\wt{i(\rho)}$ the representation of $\oJ$ given by
$$
\wt{i(\rho)}(u_-lu) = i(\rho)(l),\qquad\hbox{for }
u_-\in\oJ\cap\oU_-,\ l\in\oJ\cap\oL,\ u\in\oJ\cap\oU,
$$
whenever this defines a representation.

It will also be useful, later, to put $V'=W\oplus \{0\}\oplus W\subset V$,
equipped with the restriction of the form $\la\,,\,\ra$. We put
$A'=\End_F(V')$, $G'=\Aut_F(V')$, $\e$ the involution of
$G'$ associated to the form, and ${\oG'}^+=(G')^{\e}$. We also let
$P'$ be the maximal parabolic subgroup of $G$ which stabilizes the
flag
$$
\{0\}\subsetneq V'\subsetneq V,
$$
with unipotent radical $U'$. Let $L'$ denote the Levi
component of $P'$ which stabilizes the decomposition
$V=V'\oplus V_0$, so that $L'\simeq G'\times G_0$, and let
${P'}^{-}$ denote the opposite parabolic subgroup, with unipotent
radical ${U'}^{-}$. We note that, while $L'$ is stable under the
involution $\e$, $\oL'$ is \emph{not} a Levi subgroup of $\oG$.

\medskip

We consider the inertial class $\fs_{\oL}=
[\oL,i(\pi
)]_{\oL}$ and a type $(J_{\oL},\l_{\oL})$ for it, where
$$
J_{\oL}:= i(J_W)=i(J_W\times \oG_0)\quad\hbox{ and }\quad
\l_{\oL}:=i(\l_W
).
$$
We are going to construct a $\oG$-cover of this type, which will give a
$\oG$-type for the inertial class $\fs=[\oL,i(\pi
)]_{\oG}$.


\subsection{Semisimple characters}\label{S1.3}

We continue with the notation of the previous sections.
Let $\cL_W=\{L^W_k:k\in\bZ\}$ be the $\oe$-lattice chain in $W$ 
corresponding to $\A_W$, so that
$$
\A_W=\{a\in A_W : a L^W_k \subset L^W_k \hbox{ for all }k\in\bZ\},
$$
normalized so that $L^W_0=\of \bw_1\oplus\cdots\oplus\of \bw_M$. Let 
$e_W$ denote the $\of$-period of $\A_W$. We define $\cL'$ to be 
the $\oe$-lattice chain of $\of$-period $2e_W$ in 
$V'$ given by 
$$
\cdots\supset L^W_k\oplus L^W_k \supset L^W_k\oplus L^W_{k+1} \supset 
L^W_{k+1}\oplus L^W_{k+1} \supset\cdots
$$
 It is 
straightforward to check, since $\s\in U(\A_W)$, we have
$$
\left(L^W_k\oplus L^W_k\right)^\# = L^W_{e_W-k}\oplus L^W_{e_W-k} 
\quad\hbox{ and }\quad 
\left(L^W_k\oplus L^W_{k+1}\right)^\# = L^W_{e_W-k-1}\oplus L^W_{e_W-k} 
$$
so that $\cL'$ is a self-dual lattice chain. We write 
$\cL'=\{L'_k:k\in\bZ\}$, where we number the lattices so that 
$L'_0=L^W_{\lfloor\frac{e_W}2\rfloor}\oplus 
L^W_{\lceil\frac{e_W}2\rceil} = \left(L'_0\right)^\#$.

Now let $\L'$ be the $\oe$-lattice \emph{sequence} of period 
$4e_W$ in $V'=W\oplus W$ given by
$$
\L'(k)=L'_{\lfloor\frac k2\rfloor}
$$
so that every lattice of $\cL'$ occurs twice in the sequence 
and $\L'(k)^\#=\L'(1-k)$, for $k\in\bZ$. 

We consider the element $\frac 12\b\oplus\frac 12\b\in A_W\oplus A_W$, 
which we will also call $\frac 12\b$.  
Then $[\cL',2n_W,0,\frac 12\b]$ and $[\L',4n_W,0,\frac 
12\b]$ are both simple strata in $A'$ so we can define the 
orders $\H,\J$ and the groups $H^k,J^k$ for them, and also simple 
characters (see~\cite{BK2}). These are in fact the same, up to a scaling of 
the index ({\it loc.\ cit.}); so, for example,
$$
H^k(\tfrac 12\b,\L') = H^{\lceil\frac k2\rceil}(\tfrac 
12\b,\cL').
$$
In particular, $H^1(\tfrac 12\b,\L') = H^1(\tfrac 12\b,\cL')$ 
and we shall denote this group ${H'}^{1}$. Similarly, we put 
${J'}^{1}=J^1(\tfrac 12\b,\L') = J^1(\tfrac 12\b,\cL')$ and 
$J'=J(\tfrac 12\b,\L') = J(\tfrac 12\b,\cL')$. Likewise, 
the simple characters of ${H'}^{1}$ are the same: 
$\cC(\L',0,\frac 12\b)=\cC(\cL',0,\frac 12\b)$.

Moreover, the groups ${H'}^{1}$, ${J'}^{1}$ and $J'$ are described 
in~\cite[2.2 Lemma]{B2}. So, for example, with respect to the basis $\cB\cup\cB$ of $V'$ (which, we note, is \emph{not} a Witt basis), 
$$
{H'}^{1} = \begin{pmatrix} 
H^1_W & \J_W \\ \pie\J_W & H^1_W
\end{pmatrix}.
$$
Finally, we have the bijection $\tau_{\A_W,\A',\frac 12\b}: 
\cC(\A_W,0,\frac 12\b)\to\cC(\A',0,\frac 12\b)$ from~\cite[\S3.6]{BK}. 
Let $\th'$ be the image of $\th_L$ under this map. Note that, 
by~\cite[\S7.1--2]{BK}, $\th'$ is trivial on ${H'}^{1}\cap U$ and 
${H'}^{1}\cap U_-$, while
$$
\th'|_{{H'}^{1}\cap L}=\th_L\otimes\th_L.
$$
It is 
straightforward to check that Proposition~\ref{1.2}, together with our 
definition of the form $\la \,,\,\ra$, implies that $\th'$ 
is fixed by the involution $\e$ (cf.~\cite[2.3 Corollary]{B2}).

\medskip

Now let $\cL_0$ be the unique self-dual $\of$-lattice chain in $V_0$
$$
\cdots\supseteq L_0\supseteq L_0^\# \supseteq \pf L_0 \supseteq \pf 
L_0^\# \supseteq\cdots.
$$
Note that we may have $L_0=L_0^\#$ or $L_0^\#=\pf L_0$ so this lattice 
chain has $\of$-period $e_0=1$ or $2$. Let $\L_0$ be the self-dual 
$\of$-lattice sequence of period $4e_W$ given by
$$
\L_0(k)=\begin{cases} 
\pf^j L_0&\quad\hbox{if }\lceil\tfrac k{2e_W}\rceil = 2j, \\
\pf^j L_0^\#&\quad\hbox{if }\lceil\tfrac k{2e_W}\rceil = 2j+1,
\end{cases}
$$
so that every lattice of $\cL_0$ appears with equal multiplicity 
$4e_W/e_0$. Note also that $\L_0(k)^\#=\L_0(1-k)$, for $k\in\bZ$. 
Moreover, the filtration of $A_0$ determined by $\L_0$ is the same as 
that determined by $\cL_0$ up to a scaling of the index. In 
particular, $\fa_0(\L_0)=\A(\cL_0)$ and $\fa_1(\L_0)=\P(\cL_0)$. 

\medskip

Finally, we define $\L$ to be the $\of$-lattice sequence in $V$ given 
by
$$
\L(k)=\L'(k)\oplus\L_0(k),\qquad\hbox{for }k\in\bZ.
$$
Then, by construction, $\L$ is self-dual and of $\of$-period $e=4e_W$.
We consider the element $\frac 12\b\oplus 0 \in A'\oplus A_0$ (in 
fact, in $A_L$); by abuse of notation, we will still call this 
element $\frac 12\b$. Then $[\L,n,0,\frac 12\b]$ is a \emph{semisimple 
stratum} in $A$, where $n=4n_W$. (See~\cite[\S3.1]{S4} for the definition 
of semisimple stratum, which is more general than the definition in~\cite[\S3.3]{S1}; in particular, null strata are thought of as simple 
strata so, alternatively, the definition in~\cite{S1} could be used with 
``simple'' replaced by ``simple or null'' everywhere. The results of~\cite{S1} all remain valid in this situation -- the proofs are the same 
and they are also proved in~\cite{S4}.)

We put $J=J(\frac 12\b,\L)$ (see~\cite{S1} or~\cite{S4}) and similarly for 
$J^1$, $H^1$, $\J$, $\H$ etc. In matrix form we have
$$
\J = \begin{pmatrix} 
\J_W & \fa_{\lceil n/2\rceil}(\L) & \pie^{-1}\H^1_W \\ 
\fa_{\lceil n/2\rceil}(\L) & \A(\cL_0) & \fa_{\lceil n/2\rceil}(\L) \\ 
\H^1_W & \fa_{\lceil n/2\rceil}(\L) & \J_W
\end{pmatrix}
$$
and there are similar decompositions for $\J^k$ and $\H^k$, $k\ge
0$. Since $J'$ is stable under the involution $\e$ and $\L$ is 
self-dual, we see that $J$ is stable under the involution $\e$, and 
likewise for $J^1$ and $H^1$.

\medskip

Let $\th$ be the unique semisimple character (see~\cite[\S3.3]{S1} or~\cite[\S3.2]{S4}) of $H^1$ such that
$$
\th|_{{H'}^{1}} = \th'.
$$
Now $\th$ is trivial on $H^1\cap U'$ and $H^1\cap {U'}^{-}$ by 
definition so, since $\th'$ is trivial on ${H'}^{1}\cap U$ and 
${H'}^{1}\cap U_-$, we see that $\th$ is in fact trivial on $H^1\cap 
U$ and $H^1\cap U_-$. Since $\th|_{{H'}^{1}}$ is fixed by $\e$ and 
$\th|_{U^1(\cL_0)}$ is trivial (by definition), we see that $\th$ is 
fixed by $\e$. Moreover, since $\th_L^2=\th_W$, we have
$$
\th|_{\oH^1} = \wt{i(\th_W)}.
$$

\begin{Proposition}[{cf.~\cite[2.3 Theorem]{B2}}]\label{1.4} 
There exists a
semisimple stratum $[\L,n,0,\a]$ in $A$ with $\a\in A_{L}$ and
$\a=-\wt\e(\a)$ such that $\th\in\cC(\L,0,\a)$.
\end{Proposition}

\begin{proof} Let $\vph:A'\to A'$ be the
involution given by conjugation by 
$$
h=\begin{pmatrix} I_M&0\\0&-I_M\end{pmatrix}.
$$ 
Then $(\theta')^\varphi=\theta'$, since it is
trivial on the unipotent parts, and $(\theta')^\e=\theta'$.
Thus, $\theta'$ is invariant under the subgroup $\Omega$ of
$\Aut(G')$ generated by $\e$ and $\varphi$.
As $\e\varphi=\varphi\e$, we see $\Omega\simeq\mathbb Z_2\times\mathbb Z_2$,
and thus by~\cite[Theorem~6.3]{S1} (see the Remarks at the bottom of page
139 there also), there is a choice of $\alpha$ which is fixed by
$\Omega$ such that $\th'\in\cC(\L',0,\a)$. Then $[\L,n,0,\a]$
is a semisimple stratum with $\a\in A'\cap A_{L}$ and
$\a=-\wt\e(\a)$, such that $\th\in\cC(\L,0,\a)$. 
\end{proof}

In particular, replacing $\frac 12\b$ by the element $\a$ found in 
Proposition~\ref{1.4}, we may (and do) assume that $[\L,n,0,\frac 12\b]$ 
is a skew semisimple stratum.

\medskip

As in the construction of the type $(J_W,\l_W)$, we will also need 
the transfer of $\th$ to a minimal order, which we give now.
Let $B'$ denote the centralizer in $A'$ of $E$ and put 
$\B'=\A'\cap B'$, an $\oe$-order of period $2$ with 
radical $\Q'$. Recall that we had, in \S1.1, a minimal $\oe$-order 
$\B^\sm_W\subset\B_W$. Let ${\B'}^\sm\subset\B'$ be the 
$\wt\e$-stable minimal $\oe$-order in $B'$ given by
$$
{\B'}^\sm=(\B^\sm_W\oplus\wt\e(\B^\sm_W)) + \Q',
$$
(cf.~\cite[4.3]{B2}). 
Let ${\cL'}^\sm$ be the corresponding self-dual $\oe$-lattice chain in 
$V'$, of $\oe$-period $2r$, where $r=M/[E:F]$. Let ${\L'}^\sm$ be 
the self-dual $\oe$-lattice sequence in $V'$ in which every 
lattice of ${\cL'}^\sm$ occurs twice and with the indexing chosen 
such that
$$
{\L'}^\sm(k)^\# = {\L'}^\sm(1-k),\quad\hbox{for all }k\in\bZ.
$$
Note that $[{\cL'}^\sm,2n^{\sm}_W,0,\frac 12\b]$ and 
$[{\L'}^\sm,4n^{\sm}_W,0,\frac 12\b]$ are simple strata in $V'$ whose 
associated groups and characters are the same up to a scaling of index.

Now let $\L^\sm$ be the $\of$-lattice sequence in $V$ defined by
$$
\L^\sm(k)={\L'}^\sm(k)\oplus\L_0\left(\lceil\tfrac kr\rceil\right).
$$
It is a self-dual lattice sequence of period $re$ such that 
$\fa_0(\L^\sm)\subset\fa_0(\L)$. Then $[\L^\sm,n^{\sm},0,\frac 12\b]$ is a 
semisimple stratum in $A$, where $n^{\sm}=4n^{\sm}_W=nr$.

\medskip

We put ${\th'}^\sm=\tau_{\L',{\L'}^\sm,\frac 12\b}(\th')$, a 
simple character of $H^1(\frac 12\b,{\L'}^\sm)$; then also
${\th'}^\sm=\tau_{\A^\sm_W,{\L'}^\sm,\frac 12\b}(\th^\sm_L)$. 
Finally, let $\th^\sm$ be the unique semisimple character of 
$H^1_\sm=H^1(\frac 12\b,\L^\sm)$ such that
$$
\th^\sm|_{H^1(\frac 12\b,{\L'}^\sm)} = {\th'}^\sm.
$$
In the language of~\cite[\S3.5]{S4}, we have $\th^\sm= 
\tau_{\L,\L^\sm,\frac 12\b}(\th)$.


\Section{Covers and Hecke Algebras}\label{S2}

We continue with the notation above, so we have $[\L,n,0,\frac
12\beta]$ a skew semisimple stratum in $A$, with $\beta\in A_L$, and
$\theta\in\cC(\L,0,\frac 12\b)$ such that 
$$\theta|_{\oH^1} = \wt{i(\th_W
)}.
$$
We also have $E=F[\b]$ and $B$ the $A$-centralizer of $E$. The galois 
involution $\mu$ extends to $E$ (as $\wt\e$), and we write $E_0$ for 
the fixed subfield, which has index $2$. We fix a uniformizer $\pie$ of 
$E$ such that $\pie^\mu=\pm\pie$.


\subsection{$\boG$-Covers}\label{S2.1}

We begin by doing some work in $G$, before using Glauberman's
correspondence to transfer this information to $\oG$. 

Recall that, given $\rho$ a representation of a subgroup $H$ of $G$
and $g\in G$, the \emph{intertwining space} $I_g(\rho|H)$ is
$$
I_g(\rho|H)= \Hom_{H\cap {}^gH}(\rho,{}^g\rho),
$$
where ${}^gH=gHg^{-1}$ and ${}^g\rho$ is the representation of ${}^gH$
given by ${}^g\rho(ghg^{-1})=\rho(h)$, for $h\in H$. The 
\emph{$G$-intertwining} $I_G(\rho|H)$ of $\rho$ is then defined to be
$$
I_G(\rho|H)=\{g\in G:I_g(\rho|H)\ne 0\}.
$$

In our situation, we recall the following, from~\cite[Theorem~3.14,
Corollary~4.2 and Proposition~4.3]{S1}:

\begin{Lemma}\label{2.1}
\begin{enumerate} 
\item $I_{G}(\theta)=J^1 B^\times J^1$.
\item There exists a unique irreducible representation $\eta$ of $J^1$
which contains $\theta$. Moreover, $\dim\eta=(J^1:H^1)^{\frac12}$ and,
for $g\in G$,
$$
\dim I_g(\eta|J^1)=\begin{cases} 1 &\quad\hbox{if }g\in J^1 B^\times J^1; \\
0 &\quad\hbox{otherwise.} \end{cases}
$$
\end{enumerate}
\end{Lemma}

From~\cite[\S7.1]{BK} (together with~\cite[\S3.3]{S1}), $H^1$ has an Iwahori
decomposition with respect to $(L,P)$ and
$$
H^1\cap L = H^1(\tfrac 12\b,\A_W) \times
U^1(\L_0) \times H^1(\tfrac 12\b,\A_W).
$$
There are also similar decompositions for $J^1$ and for
$J$. Moreover, $\th$ is trivial on $H^1\cap U$
and $H^1\cap U_-$, and the restriction of $\th$ to $H^1\cap L$
takes the form 
$$
\th|_{H^1\cap L} = \th_L \otimes \bI \otimes \th_L.
$$

\medskip

We recall (\cite[Proposition~4.1]{S1} -- see also~\cite[\S3.4]{BK}) that the pairing
$$
(j,j')\mapsto \th[j,j'],\qquad\hbox{for }j,j'\in J^1
$$
induces a nondegenerate alternating bilinear form
$\bk_\th$ on $J^1/H^1$; likewise, we have a nondegenerate
alternating bilinear form $\bk_{\th_L}$ on
$J^1(\frac 12\b,\A_W)/H^1(\frac 12\b,\A_W)$. Then, exactly as in~\cite[Proposition~7.2.3]{BK}, we get:

\begin{Lemma}[{cf.~\cite[Proposition~7.2.3]{BK}}]\label{2.2} 
\begin{enumerate}
\item The subspaces
$J^1\cap U_-/H^1\cap U_-$ and 
$J^1\cap U/H^1\cap U$ of $J^1/H^1$ are both totally isotropic for the
form $\bk_\th$ and orthogonal to the subspace $J^1\cap L/H^1\cap L$.
\item The restriction of $\bk_\th$ to the group
$$
J^1\cap L/H^1\cap L = J^1(\tfrac 12\b,\A_W)/H^1(\tfrac 12\b,\A_W) 
\times J^1(\tfrac 12\b,\A_W)/H^1(\tfrac 12\b,\A_W)
$$
is the orthogonal sum of the pairings $\bk_{\th_L}$, $\bk_{\th_L}$.
\item We have an orthogonal sum decomposition
$$
\frac{J^1}{H^1} = \frac{J^1\cap L}{H^1\cap L} \perp \left(
\frac{J^1\cap U_-}{H^1\cap U_-}\times 
\frac{J^1\cap U}{H^1\cap U}\right).
$$
\end{enumerate}
\end{Lemma}

We define the groups
$$
H^1_P = H^1(J^1\cap U),\qquad J^1_P=H^1(J^1\cap P),\qquad 
J_P=H^1(J\cap P).
$$
Since $J^1$ normalizes $\th$ and $\th|_{H^1\cap U}$ is
trivial, we can define the character $\th_P$ of $H_P^1$ by
$$
\th_P(hu) = \th(h),\qquad\hbox{for }h\in H^1, u\in J^1\cap U.
$$

As in~\cite[\S7.2]{BK}, we immediately get:

\begin{Corollary}[{cf.~\cite[Propositions~7.2.4,~7.2.9]{BK}}]\label{2.3} 
There is a
unique irreducible representation $\eta_{P}$ of $J^1_{P}$ such that
$\eta_{P}|_{H^1_{P}}$ contains $\th_{P}$. Moreover,
$\eta\simeq\Ind_{J^1_P}^{J^1}\eta_{P}$ and, for each $b\in B^\times$,
there is a unique $(J^1_{P},J^1_{P})$-double coset in $J^1bJ^1$
which intertwines $\eta_{P}$.
\end{Corollary}

We note also that we certainly have $\eta_P|_{J^1_P\cap L}\simeq
\eta_L\otimes\bI\otimes\eta_L$.

\begin{Proposition}\label{2.4} We have $I_G(\th_P)= J^1_P B^\times J^1_P$ and
hence $I_G(\eta_P)= J^1_P B^\times J^1_P$.
\end{Proposition}

\begin{proof} We have $I_G(\th_P)=I_G(\eta_P)$ so, by Corollary~\ref{2.3}, we
need only check that all of $B^\times$ intertwines $\th_P$. Since
$B^\times\subset L'$, and $\th_P$ is trivial on $H^1\cap U'$ and
$H^1\cap {U'}^{-}$, we need only check that $B^\times$ intertwines
$\th_P|_{H^1\cap L'}=\th'_P\otimes\bI$, where $\th'_P$ is
the character of ${H'}^{1}_P={H'}^{1}({J'}^{1}\cap U)$ obtained by
trivial extension, as above. Since $G_0$ clearly intertwines the
trivial representation $\bI$, we need to check that ${B'}^{\times}$,
the $G'$-centralizer of $\b$, intertwines $\th'_P$. 

To ease notation, we omit the superscripts ${}'$; indeed we have
just reduced to the case where $V_0=0$. 
Then we are in the situation of~\cite[\S7.1--2]{BK}
. We note that $J_P=(U(\B)\cap L)J^1_P$. Since
$U(\B)\cap L$ normalizes $\th_L\otimes\th_L$, while $\th_P$ is trivial
on $H^1\cap U$ and $H^1\cap U_-$, we deduce that $J_P$ normalizes
$\th_P$. 

Recall that we have our ($\wt\e$-stable) minimal $\oe$-order 
$\B^\sm\subset\B$ from \S~\ref{S1.3} (where we are still omitting the 
superscripts ${}'$). 
Let $\bv_1,\ldots,\bv_{r}$ 
be an $E$-basis for $W$ and $\bv_{r+i}=\s\bv_i$, for $i=1,2,\ldots,r$, 
so that, with respect to 
the basis $\bv_1,\ldots,\bv_{2r}$ of $W\oplus W$, we have $\B^\sm$ in 
standard form. We put 
$V^{(i)}=Ev_i$, for $i=1,\ldots,2r$. Let $P_0\subset P$ denote the parabolic 
subgroup of $G$ which is the stabilizer of the flag
$$
0 \subset V^{(1)} \subset V^{(1)}\oplus V^{(2)} \subset \cdots \subset
\bigoplus_{i=1}^{2r} V^{(i)} = V 
$$
Let $U_0\supset U$ be the unipotent radical of $P_0$ and let 
$L_0\subset L$ be the Levi component of $P_0$ which stabilizes the 
decomposition
$$
V=\bigoplus_{i=1}^{2r} V^{(i)}.
$$
Also, let $P_0^-=L_0U_0^-$ denote the opposite parabolic subgroup.

For $i=1,\ldots,2r$, let $\L^{(i)}$ denote the lattice \emph{chain}
(i.e.\ ignore repetitions) given by intersection of $\L$ with
$V^{(i)}$; note that $[\L^{(i)},n,0,\frac 12\b]$ is a simple stratum in
$A^{(i)}=\End_F(V^{(i)})$. Moreover, since each $V^{(i)}$ is a
one-dimensional $E$-vector space, the lattice chain $\L^{(i)}$ is
uniquely determined by the property of being normalized by
$E^\times$. Hence, we can (and do) identify $\L^{(i)}$ with $\L^{(1)}$,
for each $i$.

From~\cite[Examples~10.9,~10.10]{BH}, $H^1$ has an Iwahori decomposition with
respect to $(L_0,P_0)$ and 
$$
H^1_P\cap L_0=H^1\cap L_0 = \prod_{i=1}^{2r} H^1(\b,\L^{(1)}).
$$
There is also a similar decomposition for $J^1$ (though not for
$J$) and hence also for $H^1_P$. Moreover ({\it loc.\ cit.}), $\th$
is trivial on $H^1\cap U_0$ and $H^1\cap U_0^-$, and the restriction
of $\th$ to $H^1\cap L_0$ takes the form 
$$
\th_P|_{H^1_P\cap L_0} = \th|_{H^1\cap L_0} = \bigotimes_{i=1}^{2r} \th^{(1)}
$$
where $\th^{(1)}$ is the simple character transfer of $\th_L$ to
$H^1(\b,\L^{(1)})$.

Since $U(\B^\sm)\subset (U(\B)\cap P)U^1(\B)$, we need only show that
some set of double coset representatives for $U(\B^\sm)\backslash
B^\times\slash U(\B^\sm)$ intertwines $\th_P$. By the Bruhat
decomposition, we may take these representatives to be of the form
$$
b=yz
$$
where, with respect to the basis $v_1,\ldots,v_{2r}$, $y$ is a
permutation matrix and $z=\diag(z_1,\ldots,z_{2r})$ is a diagonal matrix
with entries $z_i\in\pie^{\bZ}$. Then $b$ normalizes $H^1_P\cap L_0$
and indeed, since $\pie^\bZ$ normalizes $\th^{(1)}$, normalizes
$\th_P|_{H^1_P\cap L_0}$. 

Since $b$ normalizes $L_0$, $H^1_P$ also has an Iwahori decomposition
with respect to $(L_0,P_0^b)$, whence, by uniqueness of Iwahori
decompositions, $H^1_P\cap {}^bH^1_P$ has an Iwahori decomposition
with respect to $(L_0,P_0)$. Writing $g\in H^1_P\cap {}^bH^1_P$ in
this decomposition as $g=u_-lu$, we have
$$
\th_P(g)=\th_P(l)={}^b\th_P(l).
$$
To conclude the Proposition, we need only show that $\th_P(b^{-1}ub)=1$, and
likewise for $u_-$. Since both are similar, we treat only the first of
these. Since $b^{-1}ub\in H_P^1$, we can write $b^{-1}ub=u'_-l'u'$ in
the Iwahori decomposition of $H_P^1$. However, by elementary row and
column operations, we can write $b^{-1}ub=x_-x$, with $x_-\in U_-^-$
and $x\in U_0$. Now uniqueness of Iwahori decompositions implies that
$u'_-=x_-$, $l'=1$ and $u'=x$. Hence $\th_P(b^{-1}ub)=
\th_P(u'_-)\th_P(u) =1$ as required, since $\th_P$ is trivial on
$H^1_P\cap U_0$ and $H^1_P\cap U_0^-$.
\end{proof}

We will also need similar results for our character $\th^\sm$ 
of $H^1_\sm=H^1(\frac 12\beta,\L^\sm)$.
Let $\eta^\sm$ be the unique irreducible representation of
$J^1_\sm=J^1(\frac 12\beta,\L^\sm)$ which contains $\th^\sm$.
As above, $H^1_\sm$ and $J^1_\sm$ have Iwahori
decompositions with respect to $(L,P)$ and we can define the character
$\th^\sm_P$ of $H^1_{\sm,P}=H^1_\sm(J^1_\sm\cap U)$ by trivial extension of
$\th^\sm$. The same proofs (indeed, they are somewhat simpler) as those of
Lemma~\ref{2.2}, Corollary~\ref{2.3} and Proposition~\ref{2.4} show that
$$
I_G(\th^\sm_P)=J^1_{\sm,P}B^\times J^1_{\sm,P},
$$
where $J^1_{\sm,P}=H^1_\sm(J^1_\sm\cap P)$, that there is a unique
irreducible representation $\eta^\sm_P$ of $J^1_{\sm,P}$ which contains
$\th^\sm_P$, and that
$\eta^\sm=\Ind_{J^1_{\sm,P}}^{J^1_\sm}\eta^\sm_P$. 

Let $\th^\sm_L$ denote the transfer to $H^1(\frac 12\b,\A_W^\sm)$ of
$\th_L$ (this is just the restriction of $\th^\sm$ to one of the 
copies of $G_W\subset L$). Then, if $\eta^\sm_L$ denotes the unique
irreducible representation of $J^1(\frac 12\b,\A^\sm_W)$ which
contains $\th^\sm_L$, we have
$$
\eta_P^\sm|_{H^1_\sm\cap L}\simeq\eta^\sm_L\otimes\bI\otimes\eta^\sm_L.
$$


\subsection{$\ov{\boG}$-Covers}\label{S2.2}

Now we will transfer the information obtained in the last section to
$\oG$, using Glauberman's correspondence (see~\cite{G}, or~\cite{S1} for the
situation here). Let $\Om$ denote a $2$-group of automorphisms of
$G$. Recall that if $H$ is a pro-$p$ subgroup of $G$ and
$H^\Om$ is the group of $\Om$-fixed points, then there is a bijection,
denoted $\rho\leftrightarrow\bg_\Om(\rho)$ between (equivalence
classes of) irreducible representations of $H$ with
$\rho^\om\simeq\rho$, for all $\om\in\Om$, and (equivalence classes
of) irreducible representations of $H^\Om$. Further, this
correspondence commutes with irreducible restriction and irreducible
induction. Recall also that the representation
$\bg_\Om(\rho)$ is characterized as the unique component of
$\rho|_{H^\Om}$ appearing with odd multiplicity. 

We will usually apply this correspondence with $\Om$ the group of
automorphisms of $G$ consisting of $\e$ and the identity, in which
case we will just write $\bg$ for the correspondence. Note also that, for 
$\e$-stable pro-$p$ subgroups $H$ of $G$, we have $H^\e=\oH^+=\oH$.

\medskip

We write $\ov\th{}=\bg(\th)=\th|_{\oH^1}$, a skew semisimple character
(see~\cite[\S3.4]{S1}). We also set 
$$
\oG_E^+=B\cap\oG^+\hbox{ and }\oG_E=B\cap\oG
$$ 
so that $\oG_E^+$ (respectively $\oG_E$) is the direct 
product of a unitary group, for the quadratic extension $E/E_0$, and 
the anisotropic group $\oG_0^+$ (respectively $\oG_0$).

We recall first the following, 
from~\cite[Theorem~3.16]{S1} (see also the Remarks following {\it op.\ cit.} Corollary~4.2), and~\cite[Proposition~3.3.1]{S4}:

\begin{Lemma}\label{2.5}
\begin{enumerate}
\item $I_{\oG^+}(\ov\th{})=\oJ^1 \oG_E^+ \oJ^1$ and 
$I_{\oG}(\ov\th{})=\oJ^1 \oG_E \oJ^1$.
\item There exists a unique irreducible representation $\ov\eta{}$ of
$\oJ^1$ which contains $\ov\th{}$. Moreover, $\ov\eta{}=\bg(\eta)$,
$\dim\ov\eta{}=(\oJ^1:\oH^1)^{\frac12}$ and, for $g\in\oG^+$,
$$
\dim I_g(\ov\eta{}|\oJ^1)=
\begin{cases} 1&\quad\hbox{if }g\in\oJ^1\oG_E^+\oJ^1; \\
0 &\quad\hbox{otherwise.} \end{cases}
$$
\end{enumerate}
\end{Lemma}

We also put $\ov\th{}_P=\bg(\th_P)=\th_P|_{\oH^1_P}$. Note that we have 
$$
\ov\th{}_P=\wt{i(\th_W)}.
$$
We remark here that, since $\oJ_P\cap\oU=\oH^1_P\cap\oU$, and likewise
for $\oU_-$, this implies (see~\cite[Lemma~1(ii)]{B1}) that
$$
(\oJ_P\cap\oU)(\oJ_P\cap\oU_-) \subset 
(\oJ_P\cap\oU_-)\Ker(i(\th_W))(\oJ_P\cap\oU). 
$$
In particular (by {\it loc.\ cit.}), if $\rho$ is a representation of
a subgroup $K$ of $J_W\times\oG_0$ which restricts to a multiple
of $\th_W$, then $\wt{i(\rho)}$ is a well-defined representation of
$(\oJ_P\cap\oU_-)(i(K))(\oJ_P\cap\oU)$. 

\medskip

From Proposition~\ref{2.4}, together with~\cite[Corollary~2.5]{S1} and\cite[Theorem~2.3]{S2}, we get: 

\begin{Proposition}\label{2.6} $I_{\oG^+}(\ov\th{}_P)= \oJ^1_P
\oG_E^+ \oJ^1_P$ and $I_{\oG}(\ov\th{}_P)= \oJ^1_P
\oG_E \oJ^1_P$.
\end{Proposition} 

Let $\ov\th{}^\sm=\bg(\th^\sm)=\th^\sm|_{\oH^1_\sm}$. Let $\ov\eta{}$
be the unique irreducible representation of $\oJ^1$ which contains
$\ov\th{}$, as in Lemma~\ref{2.5}, and let $\ov\eta{}^\sm$ be the unique
irreducible representation of $\oJ^1_\sm$ which contains $\ov\th{}^\sm$. Let
$\ov\th{}_P=\bg(\th_P)$ be the trivial extension of $\ov\th{}$ to
$\oH^1_P$ and, similarly, $\ov\th{}^\sm_{P}$ the trivial extension 
of $\ov\th{}^\sm$ to $\oH^1_{\sm,P}$. 

We also put $\ov\eta{}_P=\bg(\eta_P)$ and
$\ov\eta{}^\sm_P=\bg(\eta^\sm_P)$. By taking $\e$-fixed points in
$J^1/H^1$, we can imitate Lemma~\ref{2.2} and Corollary~\ref{2.3} to show that
there is a unique irreducible representation of $\oJ^1_P$ which
contains $\ov\th_P$ and, since $\eta_P|_{\oH^1_P}$ is a multiple of
$\ov\th_P$, we see that this must be $\ov\eta_P$. Likewise,
$\ov\eta{}^\sm_P$ is the unique irreducible representation of
$\oH^1_{\sm,P}$ containing $\ov\th{}^\sm_P$. 
We also have
$$
\ov\eta{}_P=\wt{i(\eta_W)}\quad\hbox{ and }\quad
\ov\eta{}^\sm_P=\wt{i(\eta^\sm_W)},
$$
since the representations on the right restrict to multiples of
$\wt{i(\th_W)}=\ov\th_P$ and $\wt{i(\th^\sm_W)}=\ov\th{}^\sm_P$ respectively.

\medskip

Since $U(\L)\cap B$ normalizes $J^1_P$ and $U(\L^\sm)$ is contained in
$U(\L)$, we can form the group $\wt
J^1_P=(\oU^1(\L^\sm)\cap\oG_E)\oJ^1_P$. We 
also recall that we have the representation $\wt\eta_W$ of $\wt
J^1_W=U^1(\B_W^\sm)J^1_W$ (see \S1.2) and observe that $\wt
J^1_P\cap\oL=i(\wt J^1_W)$.

\begin{Proposition}[{cf.~\cite[Propositions~5.1.15,~5.1.19]{BK}}] \label{2.7}
There is a
unique representation $\wt\eta_P$ of $\wt J^1_P$ such that 
\begin{enumerate}
\item $\wt\eta_P|{\oJ^1_P} = \ov\eta{}_P$;
\item $\wt\eta_P$ and $\ov\eta{}^\sm_{P}$ induce equivalent
irreducible representations of $\oU^1(\L^\sm)$.
\end{enumerate}
Moreover, $\wt\eta_P=\wt{i(\wt\eta_W)}$ and
$$
\dim I_g(\wt\eta_P) = \begin{cases} 1&\hbox{ if }g\in 
\wt J_P^1\oG_E^+\wt J_P^1;
\\
0&\hbox{ otherwise.} \end{cases}
$$
\end{Proposition}

\begin{proof} Let $\Om_1$ be the group of automorphisms of $\oG^+$
generated by conjugation by
$$
h_1=\begin{pmatrix} -I_M&&\\&I_D&\\&&-I_M\end{pmatrix}.
$$
Then $(\oG^+)^{\Om_1}=(\oG')^+\times\oG_0^+$. We will denote the Glauberman
correspondence $\bg_{\Om_1}$ on representations of pro-$p$ subgroups
of $\oG$ by $\bg_1$.
We note that, since all our representations are trivial on $U'$
and ${U'}^{-}$, they are fixed by $\Om_1$ and, moreover, $\bg_1$ is
just restriction for these representations (since their restrictions
are irreducible). 

Let $\Om_2$ be the group of automorphisms of $(\oG^+)^{\Om_1}$
generated by conjugation by
$$
h_2=\begin{pmatrix} I_M&&\\&I_D&\\&&-I_M\end{pmatrix}.
$$
Then $((\oG^+)^{\Om_1})^{^{\scriptstyle\Om_2}}=\oL$. We will denote the Glauberman
correspondence $\bg_{\Om_2}$ on representations of pro-$p$ subgroups
of $(\oG^+)^{\Om_1}$ by $\bg_2$.
As above, since all our representations are trivial on $U$
and $U^-$, their transfers to $(\oG^+)^{\Om_1}$ are fixed by $\Om_2$ and,
moreover, $\bg_2$ is again just restriction for these representations
(since their restrictions are still irreducible).

We write $\bg_\Om$ for the composition $\bg_2\circ\bg_1$. Then we have
$$
\bg_\Om(\ov\eta_P)=i(\eta_W)\quad\hbox{ and }\quad
\bg_\Om(\ov\eta{}^\sm_P)=i(\eta_W^\sm).
$$

We \emph{define} $\wt\eta_P$ to be $\wt{i(\wt\eta_W)}$. Then we clearly
have $\bg_\Om(\wt\eta_P)=i(\wt\eta_W)$ and
$\wt\eta_P|_{J^1_P}=\ov\eta_P$. Also, since the intertwining of
$\wt\eta_P$ is contained in that of $\ov\eta_P$, 
$$
I_{\oG^+}(\wt\eta_P)\cap \oU^1(\L^\sm) \subset 
\oJ^1_P\oG_E^+\oJ^1_P \cap \oU^1(\L^\sm) =
(\oU^1(\L^\sm)\cap\oG_E)\oJ^1_P= \wt J^1_P,
$$
so the induced representation $\Ind_{\wt J^1_P}^{\oU^1(\L^\sm)}\wt\eta_P$
is irreducible. Likewise, $\Ind_{\oJ^1_{\sm,P}}^{\oU^1(\L^\sm)}\ov\eta{}^\sm_P$
is irreducible. Then, since Glauberman's correspondence commutes with
irreducible induction, we have
\begin{eqnarray*}
\bg_{\Om}\left(\Ind_{\wt J^1_P}^{\oU^1(\L^\sm)}\wt\eta_P\right) &=&
\Ind_{i(\wt J^1_W)}^{i(U^1(\L^\sm_W))} i(\wt\eta_W) \\ &\simeq& 
\Ind_{i(\wt J^1_{\sm,W})}^{i(U^1(\L^\sm_W))} i(\wt\eta^\sm_W) =
\bg_{\Om}\left(\Ind_{\oJ^1_{\sm,P}}^{\oU^1(\L^\sm)}\ov\eta{}^\sm_P\right).
\end{eqnarray*}
Condition {(ii)} now follows as $\bg_\Om$ is injective.

Now we show that these two conditions determine $\wt\eta_P$
uniquely. For this, we need only show that $\ov\eta{}_P$ occurs in
$\Ind_{\oJ^1_{\sm,P}}^{\oU^1(\L^\sm)}\ov\eta{}^\sm_P$ with multiplicity one. We
use the Mackey formula to compute the restriction
$$
\Res_{\oJ^1_P}^{\oU^1(\L^\sm)}
\Ind_{\wt J^1_P}^{\oU^1(\L^\sm)}\wt\eta_P.
$$
If $x\in\oU^1(\L^\sm)$ intertwines $\wt\eta_P$ with $\ov\eta{}_P$, then it
intertwines $\ov\eta{}_P$ with itself so lies in $\wt J^1_P$, as
above. Thus the multiplicity of $\ov\eta{}_P$ in
$\Ind_{\oJ^1_{\sm,P}}^{\oU^1(\L^\sm)}\ov\eta{}^\sm_P$ is equal to its
multiplicity in $\wt\eta_P$, which we know to be one.

Finally, we must show that all of $\oG_E^+$ intertwines
$\wt\eta_P$. So suppose $b\in\oG_E^+$. Since $b$ intertwines
$\ov\eta{}^\sm_P$, it certainly intertwines
$$
\Ind_{\wt J^1_P}^{\oU^1(\L^\sm)}\wt\eta_P \simeq
\Ind_{\oJ^1_{\sm,P}}^{\oU^1(\L^\sm)}\ov\eta{}^\sm_P.
$$
We deduce that there exist $u,v\in \oU^1(\L^\sm)$ such that $ubv$
intertwines $\wt\eta_P$. In particular, $ubv$ intertwines $\ov\eta{}_P$ so
there exist $j_1,j_2\in\oJ^1_P$ such that
$j_1ubvj_2\in\oG_E^+$. Note that this element also still
intertwines $\wt\eta_P$. 

Now
\begin{eqnarray*}
U^1(\L^\sm)bU^1(\L^\sm)\cap B^\times 
&=& U^1(\L^\sm)bU^1(\L^\sm)\cap L'\cap B^\times \\
&=& (U^1(\L^\sm)\cap L') b (U^1(\L^\sm)\cap L') \cap B^\times. 
\end{eqnarray*}
Now we apply~\cite[Lemma~2.1]{S2} applied with $\Gamma=\Om_1$, and we see
$$
U^1(\L^\sm)bU^1(\L^\sm)\cap B^\times = (U^1(\L^\sm)\cap B^\times)b (U^1(\L^\sm)\cap B^\times),
$$
by~\cite[Theorem~1.6.1]{BK}.

Then, applying~\cite[Lemma~2.1]{S2} again, with $\Gamma=\{1,\e\}$, we get
$$
\oU^1(\L^\sm)b\oU^1(\L^\sm)\cap\oG_E^+ = (\oU^1(\L^\sm)\cap\oG_E)b 
(\oU^1(\L^\sm)\cap\oG_E).
$$
So there exist $u',v'\in \oU^1(\L^\sm)\cap\oG_E$ such that
$u'bv'=j_1ubvj_2$ and, since $u',v'\in\wt J^1_P$, we see that 
$b$ intertwines $\wt\eta_P$. 
\end{proof}

Now we define $\ov\k_P=\wt{i(\kappa_W
)}$, a representation of
$\oJ_P$, so that $\ov\k_P|_{\wt J^1_P}\simeq\wt\eta_P$. We also put 
$\ov\rho_P=\wt{i(\rho_W)}$, a representation of $\oJ_P$ trivial on 
$\oJ_P^1$, and
$$
\ov\l_P = \wt{i(\l_W)} = \ov\k_P\otimes \ov\rho_P.
$$
Writing matrices with respect using the basis $\cB$ of $W$ in each of 
the two copies of $W$, we set
$$
w_1=\begin{pmatrix} &&\s \\ &I_D&\\ \nu\e_W(\s)&& \end{pmatrix},
$$
where $I_D$ is the $D\times D$ identity matrix.

\begin{Proposition}\label{2.8} $I_{\oG}(\ov\l_P) = \oJ_P\,
N_{\oB^\times}(\oU(\L)\cap\oG_E)\,\oJ_P$. 
\end{Proposition}

\begin{proof} Notice first ({cf.}~\cite[4.1]{B2}) that the normalizer
$N_{\oG_E}(\oU(\L)\cap\oG_E)$ is 
$(\oU(\L)\cap\oG_E)\ov\sW(\oU(\L)\cap\oG_E)$, where
$$
\ov\sW=
\left\{ i(\pie^a,\bI_{V_0}) : a\in\bZ\right\}
\cup
\left\{ i(\pie^a,\bI_{V_0}) w_1 : a\in\bZ\right\},
$$
where $\pie$ is our fixed uniformizer of $E\subset A_W$ and
$\bI_{V_0}$ is the identity map on $V_0$. 

Since the elements of $\ov\sW$ normalize $\oJ_P\cap\oL$, while $\pie$
normalizes $\l_W$ and $\l_W=\l_W\circ\wt\s$ (see Proposition~\ref{1.2}), we
see that the elements of $\ov\sW$ normalize $i(\l_W)$. On the other
hand, they either preserve $\oU$ and $\oU_-$ or interchange them so we
see that every element of $\ov\sW$ intertwines
$\ol_P=\wt{i(\l_W)}$. Hence we have $I_{\oG}(\ov\l_P) \supset \oJ_P\,
N_{\oG_E}(\oU(\L)\cap\oG_E)\,\oJ_P$.

\medskip

The proof of the opposite containment, which is a variant of the proof
of~\cite[Proposition~5.3.2]{BK}, is inspired by~\cite[page 551]{B2};
in place of~\cite[4.2 Lemma]{B2}, we use~\cite[Proposition~1.1]{S5},
which is a slight generalization of~\cite[Proposition~4.13]{M2}. It is
almost identical to the proof of~\cite[Proposition~6.14]{S5}, except
that the definition of $\ov\kappa_P$ there is {\it a priori} slightly
different. 

Suppose $g\in\oG$ intertwines $\ov\lambda_P=\ov\kappa_P\otimes\ov\rho_P$, so 
that $g\in I_{\oG}(\ov\eta_P|\oJ^1_P)=\oJ_P\oG_E\oJ_P$, as $\ov\rho_P$ 
is trivial on $\oJ^1_P$. Thus, we may assume $g$ lies in $\oG_E$. 
Moreover $\oJ_P\cap\oG_E=\oU(\L)\cap\oG_E$ is a parahoric subgroup of 
$\oG_E$ containing the Iwahori subgroup $\oU(\L^\sm)\cap\oG_E$. 
Therefore, we may further assume $g$ is a \emph{distinguished double coset 
representative} for $\oU(\L)\cap\oG_E\backslash 
\oG_E\slash\oU(\L)\cap\oG_E$ (see~\cite[\S3]{M2} or~\cite[\S1]{S5} for this 
notion).

Since $\dim I_g(\ov\eta_P,\oJ^1_P)=1$, we can imitate the proof
of~\cite[Proposition~5.3.2]{BK} to get that any non-zero intertwining
operator in  
$I_g(\ov\lambda_P,\oJ_P)$ has the form $S\otimes T$, with 
$S\in I_g(\ov\eta_P,\oJ^1_P)$ and $T$ an endomorphism of the space of 
$\ov\rho_P$.
Now the operator $S$ also intertwines the restriction 
$\kappa|\wt J^1_P = \wt\eta_P$ so, again as in~\cite[Proposition~5.3.2]{BK}, 
it follows that $T$ belongs to $I_g(\ov\rho_P|\wt J^1_P)$. In particular, 
$g$ intertwines $\ov\rho_P|_{\wt J^1_P\cap\oG_E}$. But 
$\wt J^1_P\cap\oG_E=\oU^1(\L^\sm)\cap\oG_E$ is 
the radical of the Iwahori subgroup $\oU(\L^\sm)\cap\oG_E$ of 
$\oG_E$ contained in $\oU(\L)\cap\oG_E$. By~\cite[Proposition~1.1]{S5} 
and the Remarks that follow it, we conclude that we can assume 
that $g$ normalizes $\oU(\L)\cap\oG_E$, as required.
\end{proof}

Recall that we write $J_{\oL}=i(J_W)$ and $\l_{\oL}=i(\l_W)$.
In order to prove that $(\oJ_P,\ov\l_P)$ is a cover of 
$(J_{\oL},\l_{\oL})$, the only thing remaining is to find a 
strongly $(\oP,\oJ_P)$-positive element in the centre of $\oL$ which 
supports an invertible element of the spherical Hecke algebra 
$\cH(\oG,\ov\l_P)$. To achieve this, we look at the 
$\ov\l_P$-spherical Hecke algebras of two parahoric subgroups 
whose intersections with $\oG_E$ are non-conjugate maximal 
compact open subgroups of $\oG_E$.

\medskip

Let ${\cL'}^{(1)}$ be the self-dual $\oe$-lattice chain of $\of$-period 
$e_W$ in $V'=W\oplus W$ given by 
$$
\cdots\supset L^W_k\oplus L^W_k \supset 
L^W_{k+1}\oplus L^W_{k+1} \supset\cdots,
$$
 so that 
${\cL'}^{(1)}$ consists of every second lattice of $\cL'$. Let 
${\L'}^{(1)}$ be the self-dual $\oe$-lattice sequence in $V'$ in 
which every lattice of ${\cL'}^{(1)}$ occurs twice and with the indexing 
chosen such that
$$
{\L'}^{(1)}(k)^\# = {\L'}^{(1)}(1-k),\quad\hbox{for all }k\in\bZ.
$$
Let $\L^{(1)}$ be the $\of$-lattice sequence in $V$ defined by
$$
\L^{(1)}(k)={\L'}^{(1)}(k)\oplus\L_0(2k).
$$
It is a self-dual lattice sequence of period $2e_W=e/2$ such that 
$\fa_0(\L^{(1)})\supset\fa_0(\L)$. Put $\oK_1=\oU(\L^{(1)})$.

\smallskip

We define $\oK_2=\oU(\L^{(2)})$ by the same process, starting from the
self-dual $\oe$-lattice chain ${\cL'}^{(2)}$ in $V'$ given by
$$
\cdots\supset L^W_k\oplus L^W_{k+1} \supset 
L^W_{k+1}\oplus L^W_{k+2} \supset\cdots
$$
Then $\oU(\L)\subset\oK_1\cap\oK_2$ so, in particular, 
$\oJ_P\subset\oK_1\cap\oK_2$.

Note that the element $w_1$ lies in $\ov\sW\cap\oK_1$. We also set 
$w_2=i(\pie^{-1},\bI_{V_0})w_1\in \ov\sW\cap\oK_2$.

\begin{Lemma}\label{2.9}
\begin{enumerate}
\item $\cH(\oK_1,\ov\l_P)=\la f_1,f_{w_1}\ra$ where $f_1$ is 
supported on $\oJ_P$ and $f_{w_1}$ is supported on $\oJ_P w_1 \oJ_P$.
\item $\cH(\oK_2,\ov\l_P)=\la f_1,f_{w_2}\ra$, with $f_1$ as 
in (i) and $f_{w_2}$ supported on $\oJ_P w_2\oJ_P$.
\end{enumerate}
\end{Lemma}

\begin{proof} Both parts follow from the following consideration. For 
$i=1,2$ the $\oK_i$-intertwining of $\ov\l_P$ is given by
\begin{eqnarray*} 
I_{\oK_i}(\ov\lambda_P)
&=&(\oJ_P\,N_{\oG_E}(\oU(\L)\cap\oG_E)\,\oJ_P)\cap\oK_i\\ 
&=&\oJ_P\,\left(N_{\oG_E}(\oU(\L)\cap\oG_E)\cap\oK_i\right)\,\oJ_P.
\end{eqnarray*} 
But
$(N_{\oG_E}(\oU(\L)\cap\oG_E)\cap\oK_i=\{1,w_i\}$. Moreover,
the restriction of $\ol_P$ to $\oJ_P\cap\oL$ is irreducible and $w_i$
normalizes $\oJ_P\cap L$ so the intertwining space $I_{w_i}(\ol_P)$ is
$1$-dimensional. 
\end{proof}

\begin{Lemma}\label{2.10} Consider $f_{w_1},f_{w_2}$ as elements of 
$\cH(\oG,\ov\l_P)$. Then the convolution $f_{w_1}*f_{w_2}$ is supported on
$\oJ_P w_1 w_2 \oJ_P$.
\end{Lemma}

\begin{proof} We know $f_{w_1}*f_{w_2}$ is supported on 
\begin{eqnarray*}
&&\oJ_p w_1 \oJ_P w_2 \oJ_P \\
&&\qquad=\oJ_P 
\left(w_1(\oJ_P\cap\oU_-)w_1^{-1}\right)
\left(w_1(\oJ_P\cap\oL)w_1^{-1}\right) 
w_1 w_2 
\left(w_2^{-1}(\oJ_P\cap\oU)w_2\right) \oJ_P.
\end{eqnarray*}
Since $\oJ_P$ contains $w_1(\oJ_P\cap\oU_-)w_1^{-1}$, 
$w_1(\oJ_P\cap\oL)w_1^{-1}$ and 
$w_2^{-1}(\oJ_P\cap\oU)w_2$, the lemma follows.
\end{proof}

We will prove that $f_\z:=f_{w_1}*f_{w_2}$ is invertible. 
To accomplish this we prove that $f_{w_1}$ and $f_{w_2}$ are each
invertible. In each case we know
\begin{equation}\label{2.11}
f_{w_i}* f_{w_i}=c_{1i} f_1+d_i f_{w_i}
\end{equation}
by Lemma~\ref{2.9}. Thus, we only need to show $c_{1i}\not= 0$ for each $i$.

\begin{Lemma}\label{2.12} In equation \eqref{2.11} $c_{1i}\not= 0$ for $i=1,2$.
\end{Lemma}

\begin{proof} We treat only the case $i=1$, since the other case is 
identical. We need to check that 
$f_{w_1}*f_{w_1}(1)\ne 0$. Since 
$$
f_{w_1}(x)=\begin{cases} 0&\text{if $x\not\in\oJ_P w_1\oJ_P$}\\
\ol_P\spcheck(j_1) f_{w_1}(w_1)\ol_P\spcheck(j_2)&\text{if 
$x=j_1 w_1 j_2$, with $j_1,j_2\in\oJ_P$,}\end{cases}
$$
we can write
\begin{eqnarray*} 
f_{w_1}*f_{w_1}(1)&=
&\int_{\oK_1} f_{w_1}(y) f_{w_1} (y^{-1}) dy \\ 
&=&\frac1{\left|\oJ_P\cap {}^{w_1}\!\oJ_P\right|} 
\int_{\oJ_P\times\oJ_P} 
f_{w_1}(j_1 w_1 j_2) f_{w_1} (j_2^{-1} w_1^{-1} j_1^{-1}) dj_1 dj_2\\ 
&=&\frac{\left|\oJ_P\right|}{\left|\oJ_P\cap{}^{w_1}\!\oJ_P\right|} 
\int_{\oJ_P} \ov\l_P\spcheck(j_1)f_{w_1}(w_1)f_{w_1}(w_1^{-1}) 
\ov\l_P\spcheck(j_1^{-1})) dj_1.
\end{eqnarray*}
Now $w_1$ intertwines $\ov\l_P\spcheck=\wt{i(\l_W\spcheck)}$ and 
normalizes $\oJ_P\cap\oL=i(J_W)$ so $f_{w_1}(w_1)$ is an 
equivalence ${}^{w_1} i(\l_W\spcheck)\simeq i(\l_W\spcheck)$. But 
$w_1^{-1}=h_1^{-1}w_1$, where $h_1=i(\nu\s\e_W(\s),\bI_{V_0})$ and, by 
Proposition~\ref{1.2}, $\nu\s\e_W(\s)\in J_W$; hence $w_1^{-1}\in J_P w_1$
and $f_{w_1}(w_1^{-1})$ is an 
equivalence $i(\l_W\spcheck)\simeq {}^{w_1}i(\l_W\spcheck)$. Thus 
$f_{w_1}(w_1) f_{w_1}(w_1^{-1})$ is an equivalence of 
$i(\l_W\spcheck)$ and hence a scalar $c\ne 0$. Thus
$$
f_{w_1}*f_{w_1}(1)= c\frac{\left|\oJ_P\right|^2}
{\left|\oJ_P\cap {}^{w_1}\!\oJ_P\right|}\neq 0.
$$
Therefore, $f_{w_1}$ is invertible, as required.
\end{proof}

\begin{Lemma}\label{2.13} For each $k\in\bN$, the $k$-fold convolution 
$f_\z^k$ is supported on $\oJ_P\z^k\oJ_P$.
\end{Lemma}

\begin{proof} This is by simple induction on $k$, since
\begin{eqnarray*}
&&\oJ_p \z \oJ_P \z^k \oJ_P \\
&&\qquad=\oJ_P 
\left(\z(\oJ_P\cap\oU_-)\z^{-1}\right)
\left(\z(\oJ_P\cap\oL)\z^{-1}\right) 
\z^{k+1} 
\left(\z^{-k}(\oJ_P\cap\oU)\z^k\right) \oJ_P,
\end{eqnarray*}
while $\oJ_P$ contains $\z(\oJ_P\cap\oU_-)\z^{-1}$, 
$\z(\oJ_P\cap\oL)\z^{-1}$ and 
$\z^{-k}(\oJ_P\cap\oU)\z^k$ (cf.~Lemma~\ref{2.10}).
\end{proof}

In particular, $f_\z^{e(E/F)}$ is an invertible element of 
$\cH(\oG,\ol_P)$ which is supported on the double coset
$\oJ_P\z^{e(E/F)}\oJ_P=\oJ_P \z_F\oJ_P$, where
$$
\z_F= i(\pif\bI_W,\bI_{V_0})
$$
is a strongly $(\oP,\oJ_P)$-positive element of the centre of $\oL$. 
We conclude:

\begin{Theorem}\label{2.14} Let $\pi_W$ be an irreducible supercuspidal 
representation of $G_W\simeq GL_M(F)$, with $\pi_W^{\e_W}\simeq\pi_W$. 
Using the notation above, the pair $(\oJ_P,\ol_P)$ is a $\oG$-cover 
of $(J_{\oL},\l_{\oL})$. In particular, it is an $\fs$-type, with 
$\fs=[\oL,i(\pi_W)]_{\oG}$.
\end{Theorem}

\begin{proof} By construction, $\oJ_P$ is decomposed with respect to $\oP$ so 
(i) and (ii) of~\cite[Definition~8.1]{BK1} are satisfied for $\oP$, and 
likewise for $\oP_-$. By Lemma~\ref{2.12} the strongly $(\oP,\oJ_P)$-positive 
element $\z_F$ supports an invertible element of $\cH(\oG,\ol_P)$. 
Thus $(\oJ_P,\ol_P)$ also satisfies (iii) of~\cite[Definition~8.1]{BK1} for 
$\oP$, and is therefore a $\oG$-cover of $(J_{\oL},\l_{\oL})$. 
Then~\cite[Theorem~8.3]{BK1} implies $(\oJ_P,\ol_P)$ is an $\fs$-type.
\end{proof}


\subsection{Hecke Algebras}\label{S2.3}

In this section we derive results analogous to those of Chapter 5 of~\cite{BK}. In particular, we show that the Hecke algebra of our type,
$\cH(\oG,\ol)$ can be computed by using analogous
computations in a case where the representation $\pi$ of $\oL$ is of
level zero. For many of these situations,~\cite{KM} will give us the
parameters of the Hecke algebra. 

We fix $i=1$ or $2$ and, in the notation of the previous section, we 
put $w=w_i$.
Note that, for $g\in G_W$ and $g_0\in\oG_0$, we have ${}^{w_1} i(g,g_0)=i(\wt\s(g),g_0)$ (and a similar result for $w_2$) so,
by Proposition~\ref{1.2}, we have $^{w}i(\th_W)=i(\th_W)$,
and similarly for $\eta_W$. In particular, $w$ intertwines the 
representation $\ov\eta_P=\wt{i(\eta_W)}$ of $\oJ_P$.

Now $\k_W\circ\wt\s$ is also a $\b$-extension of
$\eta_W$ so $\k_W\circ\wt\s\simeq\k_W\otimes\chi_W$, for some character
$\chi_W$ of $U(\B_W)/U^1(\B_W)\simeq J_W/J^1_W$ which factors through
the determinant $\det_{B_W/E}$. Then 
$$
^w i(\k_W)\simeq i(\k_W\circ\wt\s)=i(\k_W)\otimes i(\chi_W)
$$
and $w$ intertwines $\ov\k_P$ with $\ov\k_P\otimes\ov\chi_P$, where
$\ov\chi_P=\wt{i(\chi_W)}$. 
Note also that, since $\oJ/\oJ^1\simeq\oJ_P/\oJ_P^1\simeq J_W/J^1_W$,
we can extend $\ov\chi_P$ to the character $\ov\chi=\wt{i(\chi_W)}$ of
$\oJ$. Set $\ov\k{}=\Ind_{\oJ_P}^{\oJ}\ov\k{}_P$;
since $I_{\oG}(\ov\k{}_P)=\oJ_P\oG_E\oJ_P$,
we know that $\ov\k{}$ is irreducible and further, 
by~\cite[Proposition~4.1.3]{BK}, $w$ intertwines $\ov\k{}$ with
$\ov\k{}\otimes\ov\chi$.

\medskip

In the notation of the previous section, we put $\L^\sM=\L^{(i)}$. Then
$[\L^{\sM},n_\sM,0,\frac 12\b]$ is a skew semisimple stratum and we set 
$\oJ_\sM=\oJ(\frac 12\b,\L^\sM)$, and similarly $\oJ^1_\sM$ and $\oH^1_\sM$.
Let $\ov\th^\sM$ denote the transfer of $\ov\th$ to $\oH^1_\sM$.
By~\cite[Corollary~4.2]{S1}, there is a unique irreducible
representation $\ov\eta^\sM$ of $\oJ^1_\sM$ which extends $\ov\th_\sM$,
and by~\cite[Theorem~4.1]{S5}, we may choose a $\b$-extension,
$\ov\k_\sM$ of $\ov\eta_\sM$ to $\oJ_\sM$ -- we recall here what we 
mean by $\b$-extension: 

Recall the lattice sequence $\L^\sm$ from~\S\ref{S1.3} such that 
$\oU(\L^\sm)\cap\oG_E$ is an Iwahori subgroup of $\oU(\L)\cap\oG_E$; then 
we have $\ov\th^\sm$ the transfer of $\ov\th$ and $\ov\eta_\sm$ the unique 
irreducible representation of $\oJ^1(\frac 12\b,\L^\sm)$ 
containing $\ov\th^\sm$. We abbreviate
$\oU(\L)\cap\oG_E=\oU(\L_{\oe})$ (and similarly for other lattice
sequences) and define $\hat J_\sM=\oU(\L_{\oe})\oJ^1_\sM$ and $\hat
J_\sM^1=\oU^1(\L^\sm_{\oe})\oJ^1_\sM$. Thus,
$\oJ_\sM\supset \hat J_\sM\supset \hat J_\sM^1 \supset\oJ_\sM^1$ and
$\hat J_\sM^1$ is a pro-$p$ Sylow subgroup of
$\oJ_\sM$. By~\cite[Proposition~3.7]{S5}, there is a 
unique irreducible representation, $\hat\eta_\sM$ 
of $\hat J_\sM^1$ which extends $\ov\eta_\sM$ and such that 
$\hat\eta_\sM$ and $\ov\eta_\sm$ induce equivalent irreducible representations
of $\oU^1(\L_\sm)$. Then a $\b$-extension $\ov\k_\sM$ is an extension to 
$\oJ_\sM$ of $\hat\eta_\sM$.

Similarly (as in~\cite[Proposition~5.2.5]{BK} -- 
see~\cite[Lemma~4.2]{S5}), there is a
unique irreducible representation, $\hat\mu_\sM$ 
of $\hat J_\sM$ which extends $\ov\eta_\sM$ and such that 
\begin{equation}\label{eq2.12}
\Ind_{\hat J_\sM}^{\oU(\L_{\oe})\oU^1{(\L)}}\hat\mu_\sM\simeq
\Ind_{\oJ}^{\oU(\L_{\oe})\oU^1{(\L)}}\ov\k{} \simeq
\Ind_{\oJ_P}^{\oU(\L_{\oe})\oU^1{(\L)}}\ov\k_P{}.
\end{equation}
Moreover, as in~\cite[Proposition 5.2.6]{BK}, we have $\hat \mu_\sM|_{\hat J_\sM^1}=\hat\eta_\sM$ and, as in\cite[Proposition~2.9]{Se}, the 
$\oG_E$-intertwining of $\hat\mu_\sM$ and $\ov\k_P$ are the same; 
indeed, the same proof shows that $w$ intertwines
$\hat\mu_\sM$ with $\hat\mu_\sM\otimes\hat\chi_\sM$, where
$\hat\chi_\sM=\wt{i(\chi_W)}$ on $\hat J_\sM$.

\begin{Lemma}\label{2.15}
There is a choice of $\k{}_W$ for which $\k_W\circ\wt\s\simeq\k_W$.
\end{Lemma}

\begin{proof}
Fix some choice of $\k_W$, which fixes the character $\chi_W$. Recall that 
we have 
$\ov\k_\sM|_{\hat J^1_\sM}=\hat\eta_\sM=\hat\mu_\sM|_{\hat J^1_\sM}$, by
construction of $\beta$-extensions in~\cite[Theorem~4.1]{S5}.
Thus, if we let $\hat\k_\sM=\ov\k_\sM|_{\hat J_\sM}$, we must have
$\hat\k_\sM\simeq\hat\mu_\sM\otimes\ov\psi_\sM$ for a character
$\ov\psi_\sM$ of $\oU(\L_{\oe})/\oU^1(\L^\sm_{\oe})$, that is, a
character of the Siegel Levi subgroup of
$\oU(\L^\sM_{\oe})/\oU^1(\L^\sM_{\oe})$ which is trivial on the  maximal unipotent radical $\oU^1(\L^\sm_{\oe})/\oU^1(\L^\sM_{\oe})$. Then 
$\ov\psi_\sM$ factorizes
through the determinant on the Levi subgroup and we can write
$\ov\psi_\sM=\wt{i(\psi_W)}$, for some character $\psi_W$ of $J_W/J^1_W$.

Now $w\in \oJ_\sM$ so $w$ certainly intertwines $\hat\k_\sM$ with
itself. Hence $w$ intertwines $\hat\mu_\sM\otimes\ov\psi_\sM$ with
$\hat\mu_\sM\otimes\hat\chi_\sM\otimes{}^{w}\ov\psi_\sM$. Chasing back through
the constructions above, we see that $w$ intertwines
$\ov\k_P\otimes\ov\psi_\sM$ with
$\ov\k_P\otimes\hat\chi_\sM\otimes{}^{w}\ov\psi_\sM$. Since $w$
normalizes $\oL$, this implies that conjugation by $w$ gives an
equivalence
$\k_W\otimes\psi_W\simeq\k_W\otimes\chi_W(\psi_W\circ\wt\s)$. By~\cite[Theorem~5.2.2]{BK}, we deduce that 
$\psi_W=\chi_W(\psi_W\circ\wt\s)$ and,
in particular, $\k_W\otimes\psi_W$ is a $\b$-extension with the required
property.
\end{proof} 

We now choose $\k_W$ as in Lemma~\ref{2.15} and make the same constructions as
before: $\ov\k_P$, $\ov\k$ and $\hat\mu_\sM$. Comparing $\hat\mu_\sM$
with $\hat\k_\sM$, we again get a character $\ov\psi_\sM$ but now with the
property that ${}^{w}\ov\psi_\sM=\ov\psi_\sM$. If $E/E_0$ is unramified then
this implies that $\ov\psi_\sM={}^{w}\ov\psi_0\ov\psi_0$, for some character
$\ov\psi_0$. Then $\ov\psi_0$ extends to a character of $J_\sM$ and, 
replacing $\k_\sM$ by $\k_\sM\otimes\ov\psi_0^{-1}$, we may assume
$\ov\psi_\sM=1$. In the ramified case, the condition on $\ov\psi_\sM$ is
that $\ov\psi_\sM^2=1$ but it is (at least in principle) possible that
$\ov\psi_\sM\ne 1$. 

We write $\ov\psi_\sM=\wt{i(\psi_W)}$; then we may assume $\psi_W^2=1$.


\begin{Proposition}\label{2.16} With the notation as above, we have a 
support-preserving algebra isomorphism
$$
\cH(\oJ_\sM,\ol_P)\simeq
\cH(\oU(\L^\sM_{\oe})/\oU^1(\L^\sM_{\oe}),\wt{i(\rho_W\otimes\psi_W)}).
$$
\end{Proposition}

\begin{proof}
From the isomorphisms in~\eqref{eq2.12} 
and~\cite[Corollary~4.1.5]{BK},
we get a support-preserving isomorphism 
$$
\cH(\oJ_\sM,\ol_P)
\simeq\cH(\oJ_\sM,\ov\k{}\otimes\wt{i(\rho_W)}).
$$
Similarly, as in~\cite[Proposition~5.5.13]{BK}, since $\psi_W^2=1$ we have a
support-preserving isomorphism
$$
\cH(\oJ_\sM,\ov\k{}\otimes\wt{i(\rho_W)})
\simeq\cH(\oJ_\sM,
(\hat\mu_\sM\otimes\ov\psi_M)\otimes(\wt{i(\rho_W\otimes\psi_W)}).
$$
Finally, we have support-preserving isomorphisms
\begin{eqnarray*}
\cH(\oJ_\sM,
(\hat\mu_\sM\otimes\ov\psi_M)\otimes(\wt{i(\rho_W\otimes\psi_W)}) 
&\simeq&\cH(\oJ_\sM,\wt{i(\rho_W\otimes\psi_W)})\\
&\simeq&\cH(\oU(\L^\sM_{\oe})/\oU^1(\L^\sM_{\oe}),\wt{i(\rho_W\otimes\psi_W)}),
\end{eqnarray*}
where the first isomorphism follows from the fact that
$\hat\mu_\sM\otimes\ov\psi_\sM=\hat\k_\sM$ extends to a representation
$\k_\sM$ of $\oJ_\sM$ (cf.~\cite[Lemma~5.6.3]{BK}), and the second by
reduction modulo $\oJ^1_\sM$, since $\oJ_\sM/\oJ^1_\sM\cong
\oU(\L^\sM_{\oe})/\oU^1(\L^\sM_{\oe})$. Putting these isomorphisms together
gives the isomorphism of the Proposition.
\end{proof}

\begin{Remarks}\rm\begin{enumerate}
\item Note that, writing $\B'_\sM$ for the self-dual $\oe$-order
$\fa_0(\L^\sM)\cap B'$, we have 
\begin{eqnarray*}
\oU(\L^\sM_{\oe})/\oU^1(\L^\sM_{\oe}) &\cong& 
\oU(\B'_\sM)/\oU^1(\B'_\sM) \times \oG_0/\oG^1_0
\quad\hbox{ and } \\
\oU(\L_{\oe})/\oU^1(\L^\sM_{\oe}) &\cong &
\oU(\B')/\oU^1(\B'_\sM) \times \oG_0/\oG^1_0,
\end{eqnarray*}
where $\oG^1_0$ is the pro-$p$ radical of the anisotropic group
$\oG_0$. Then we have an isomorphism
$$
\cH(\oU(\L^\sM_{\oe})/\oU^1(\L^\sM_{\oe}),\wt{i(\rho_W\otimes\psi_W)})\simeq
\cH(\oU(\B'_\sM)/\oU^1(\B'_\sM),\wt{i(\rho_W\otimes\psi_W)}).
$$
The quotient $\oU(\B'_\sM)/\oU^1(\B'_\sM)$ is a unitary (if
$E/E_0$ is unramified), symplectic or orthogonal group over
$k_{E_0}$ and the Hecke algebra on the right is described
in~\cite{HL}. Alternatively, reduction modulo $\oU^1(\B'_\sM)$ gives 
a support-preserving isomorphism 
$$
\cH(\oU(\B'_\sM)/\oU^1(\B'_\sM),\wt{i(\rho_W\otimes\psi_W)})\simeq
\cH(\oU(\B'_\sM),\wt{i(\rho_W\otimes\psi_W)}),
$$
and the latter is described in~\cite{M1}.
\item Since (in the case where $E/E_0$ is ramified) we have
$\chi_W^2=1$, we may replace our choice of $\k_W$ by
$\k_W\otimes\chi_W$ (which has the same property of being fixed by
$\wt\s$); this replaces $\rho_W$ by $\rho_W\otimes\chi_W$, another
self-dual cuspidal representation of $\oU(\B_\sM)$, and we lose the 
character $\chi_W$ from the RHS of the isomorphism in
Proposition~\ref{2.16}. \emph{However}, we cannot do this
independently for the two choices $\L^{(1)},\L^{(2)}$ for $\L^\sM$. In
particular, if we choose (as we always can) to dispose with the
character $\chi_W$ in one case, then it may still be non-trivial in
the other.
\item Since $(\oJ_P,\ol_P)$ is a cover of $(J_{\oL},\l_{\oL})$,
by~\cite[Corollary~7.12]{BK1} we have a canonical embedding of Hecke algebras 
$t_P\!:\!\cH(\oL,\l_{\oL})\hookrightarrow\cH(\oG,\ol_P)$ and we identify
$\cH(\oL,\l_{\oL})$ (which is just the algebra of Laurent polynomials
in a single variable) with its image $\mathcal B_P$. We also put
$\mathcal K=\cH(\oJ_\sM,\ol_P)$. 
Then~\cite[Theorem~1.5]{BKcovers} implies that the map
\begin{eqnarray*}
\mathcal B_P\otimes_{\bC} \mathcal K &\to& \cH(\oG,\ol_P) \\
f\otimes\phi &\mapsto& f*\phi
\end{eqnarray*}
is an isomorphism of $(\mathcal B_P,\mathcal K)$-bimodules.
\end{enumerate}
\end{Remarks}


\end{document}